\documentclass[12pt]{article}
\usepackage{amsmath,latexsym,amsthm}
\usepackage{amssymb}
\usepackage{amsfonts}
\usepackage{bbm} 
\newtheorem{thm}{Theorem}[section]

\newtheorem{pro}[thm]{Proposition}
\newtheorem{cor}[thm]{Corollary}

\begin{document}

\title{Injective envelopes of real C*- and AW*-algebras.}

\author{Abdugafur \textsc{Rakhimov}, \ Laylo \textsc{Ramazonova}}



\begin{center}
{\bf Injective envelopes of real C*- and AW*-algebras.}\\[2mm]

{\bf Rakhimov A.A., \ Ramazonova L.D.}\\[4mm]
\end{center}

{\small
\hspace{1.0cm} National University of Uzbekistan, Tashkent, Uzbekistan,

\hspace{1.0cm} Karadeniz technical university, Trabzon, Turkey. rakhimov@ktu.edu.tr

\hspace{1.0cm} V.I.Romanovsky Institute of Mathematics of the Academy of Sciences of Uzbekistan

\hspace{1.0cm} rlaylo2405@gmail.com
}

\vspace{1.0cm}

{\footnotesize

\hspace{1.0cm} 2020 Mathematics Subject Classification. Primary 46L35, 46L10; Secondary 47L30, 47L40 

\hspace{1.0cm} Key Words and Phrases. C*- and AW*-algebras, injective envelopes of real C*-algebras

}


\begin{abstract}
It is shown that every outer *-automorphism of a real C*-algebra can be uniquely extended to an injective envelope of real C*-algebra.
It is proven that if a real C*-algebra is a simple, then its injective envelope is also simple, and it is a real AW*-factor.
The example of a real C*-algebra that is not real AW*-algebra and
the injective envelope is a real AW*-factor of type III, which is not a real W*-algebra is constructed.
This leads to the interesting result that
up to isomorphism, the class of injective real (resp. complex) AW*-factors of type III is at least one larger than the class
injective real (resp. complex) W*-factors of type III.
\end{abstract}

\section{Introduction}
\label{Intro}

\def\e{{\bf 1}\!\!{\rm I}}

As is known, injective (complex and real) W*- and C*-algebras,
in particular, factors have been studied quite well.
On the other hand, in an arbitrary case, i.e., in the non-injective case,
it is quite difficult to study (up to isomorphism) the W*-algebras ,
in particular, it is known that there is a continuum of pairwise
non-isomorphic non-injective factors of type II$_1$. Therefore,
it seems interesting to study the so called maximal injective W*
and C*-subalgebras or what is equivalent, the smallest injective
C*-algebra containing a given algebra, which is called an injective envelope of C*-algebra.
In the complex case, such algebras were considered in the works  M.Hamana, K.Saito, M.Wright,
R.Kadison, J.Fang. In this paper, we consider the existence and uniqueness of injective envelope
real C*-algebras for given real C*-algebras.
By analogy with \cite{Saito84} we will show that every outer *-automorphism of a real
C*-algebra can be uniquely extended to an injective envelope of real C*-algebra.
Similar to \cite{Hamana79}, we will prove that if a real C*-algebra is a simple, then
its injective envelope is also simple and it is a real AW*-factor.
We will build an example of a real C*-algebra that is not real AW*-algebra and
the injective envelope is a real AW*-factor of type III, which is not a real W*-algebra.
This leads to the interesting result that
up to isomorphism, the class of of injective real (resp. complex) AW*-factors of type III is at least one larger than the class
injective real (resp. complex) W*-factors of type III.

\section{Preliminaries}

Let $B(H)$ be the algebra of all bounded linear operators acting on a complex Hilbert space $H$.
Recall that a real *-subalgebra $R\subset B(H)$ with the  identity $\e$ is called a
{\it real W*-algebra}, if it is weakly closed and $R\cap iR=\{0\}$.

Let $A$ be a Banach *-algebra over the field $\mathbb{C}$. The algebra $A$ is called a C*-{\it algebra}, if
$\|aa^*\|=\|a\|^2$, for any $a\in A$. A real Banach *-algebra $R$ is called a {\it real} C*-{\it algebra},
if $\|aa^*\|=\|a\|^2$ and an element $\e + aa^*$ is invertible for any $a\in R$. It is easy to see that $R$
is a real C*-algebra if and only if a norm on $R$ can be extended onto the complexification $A=R+iR$ of
the algebra $R$ so that algebra $A$ is a C*-algebra (see \cite{Ayupov2010}, \cite{Ayupov97}).

Denote by $M_n(A)$ algebra of all $n\times n$ matrices over $A$ which is also a C*-algebra,
relative to ordinary matrix operations. An element $a\in A$ is called {\it positive} and
denote by $a\geq 0$, if there exists a self-adjoint element $b\in A$, such that $a=b^2$.
The set of all positive elements of $A$ denoted by $A^+$.

A continuous linear  map $\varphi$ between two C*-algebras $A$ and $B$ is called {\it completely positive},
if for any $n\geq 1$, the natural map $\varphi _n$ from the C*-algebra $A\otimes M_n$ to the C*-algebra $B \otimes M_n$, defining by
$$
\varphi _n\bigl( (a_{ij})_{i,j=1}^n\bigr) = \bigl(\varphi(a_{ij})\bigr)_{i,j=1}^n
$$
is positive, where $M_n$ is the C*-algebra of $n\times n$ matrices over $\mathbb{C}$.

A C*-algebra $A$ with the unit $\e$ is said to be \textit{injective} if
whenever $B$ is a unital C*-algebra and $S$ is a self-adjoint subspace of $B$ containing the unit,
then each completely positive map $\varphi : S\to A$ can be extended to a completely positive map ${\overline{\varphi}} : B\to A$
(\cite{Loebl74}, \cite{Arveson69}). W.Arveson proved \cite{Arveson69}, that for each Hilbert space $H$,
the algebra of bounded operators on $H$ is injective, i.e., $B(H)$ is injective.\\[-3mm]

It is known that any C*-algebra can be isomorphically embedded into some algebra $B(H)$ in which it is uniformly closed.
In 1967 in the work of Hakeda-Tomiyama, C*-algebras with property E are considered:
A C*-algebra $A$ has {\it the property} E (an extensional property) if there is a projection $P : B(H)\to A$ such that $\|P\|=1$ and $P(\e)= \e$.
In this case the map $P$ is completely positive. In \cite{Loebl74} R.Loebl showed equivalence of these definitions.\\[-2mm]

Thus, every C*-algebra lies in the injective C*-algebra $B(H)$ and the smallest injective C*-algebra containing this algebra is called an {\it injective enveloping} C*-algebra.
In \cite{Hamana79}, the author proved that every C*-algebra $A$ with identity has a unique injective enveloping C*-algebra, which is denoted as $I(A)$, i.e.,
the smallest injective C*-algebra containing $A$ as a C*-subalgebra.

\section{Outer *-automorphisms of real C*-algebras}

Let $A$ be a real or complex algebra. A subspace $I$ of an algebra $A$ is called an left ideal (resp. right ideal)
if $xy \in I$ (resp. $yx\in I$), for all $x\in A$ and $y\in I$. A left and right ideal is called a two-sided ideal or ideal.
An algebra $A$ is said to be $\textit{simple}$ if it contains no non-trivial two-sided ideals and the multiplication operation is not zero (that is,
there are $a$ and $b$ with $ab\neq 0$).
If $\alpha$ is an automorphism of $A$, then $A$ said to be $\alpha$-simple if the only an $\alpha$-invariant, closed, proper, two-sided ideal of $A$ is 0.\\[-2mm]

Let $A$ be (complex) an *-algebra. A real subspace $I$ of $A$ is called a \textit{real ideal of} $A$
if \ $I\cdot A, A\cdot I \subset I_c$, where $I_c=I+iI$.
Since each complex subspace of $A$ is a real subspace, any complex ideal is automatically a real ideal of $A$.
Let $I$ be a real ideal of $A$. If there exists a real *-subalgebra $R$ of $A$ with $R+iR=A$,
such that $I \subset R$, then $I$ is called a $\textit{pure real ideal}$ of $A$.
In this case, it is obvious that we have $I\cdot R \subset I$. Note that, the reverse is not true, i.e., from $I\cdot R \subset I$ it does not follow $I \subset R$. But a complex subspace $J=I+iI$ always is a complex ideal of $A$.
On the other hand if $I\subset R$ is a real subspace of $A$ and $I+iI$ is a complex ideal,
then $I$ is a pure real ideal, i.e., we obtain $I\cdot R \subset I$.

Let, now $I$ and $Q$ be pure real ideals of $A$.
In general, the set $I+iQ$ is not a (complex) subspace. More precisely the set $I+iQ$ is a complex subspace if and only if $I=Q$.
Therefore we consider the smallest complex subspace $J$ of $A$, containing $I$ and $Q$. Obviously $J$ is equal to $(I+Q)+i(I+Q)$.
Thus, if $I$ and $Q$ are real ideals, then $J=(I+Q)+i(I+Q)$ is a complex ideal.

\begin{thm}
Let $R$ be real a algebra. Then $R$ is simple if, and only if $R+iR$ is simple.
\end{thm}

\begin{proof}
Assume that real subalgebra $I$ of a algebra $R$ be  a non-trivial two-sided ideal of a algebra $R$. Then, obviously that complex subspace $I+iI$ is be a non-trivial two-sided ideal of a algebra $R+iR$. Inversely, let $J$ is be non-trivial two-sided ideal in $R+iR$ and let $a+ib \in J$, where $a,b\in R$. Since $J$ is a two-sided ideal, then for $a-ib\in {R+iR}$ we have
$$x=(a+ib)(a-ib)=a^2+b^2+iba-iab\in J$$
$$y=(a-ib)(a+ib)=a^2+b^2-iba+iab\in J$$
Hence, $x+y=2a^2+2b^2\in J$, therefore $J\bigcap R\not=\varnothing$. Assume that $I:=J\bigcap R$. It is easy to see, that real subspace $I$ be a non-trivial two-sided ideal of algebra $R$. The proof is completed.
\end{proof}
Let $\alpha$ be a *-automorphism of real *-algebra $R$. By $\tilde{\alpha}$ we denote
the linear extension of $\alpha$ to $A=R+iR$, which is defined as $\tilde{\alpha}(x+iy)=\alpha(x)+i\alpha(y)$, where $x,y\in R$.

\begin{pro}
  \label{pro:1}
$R$ is an $\alpha$-simple if and only if $A$ is an  $\tilde{\alpha }$-simple.
\end{pro}

\begin{proof}
Let $I$ be a subspace of $R$. Then $I_c=I+iI$ is (complex) subspace of $A$.
It's obvious that $I$ is ideal of $R$ if and only if $I_c$ is ideal of $A$. Moreover,
$I$ is closed if and only if $I_c$ is closed.
It's easy to see that if $\alpha(I)\subset I$, then $\tilde {\alpha}(I_c)=\alpha(I)+i\alpha(I)\subset I_c$, and
obversely, $\tilde{\alpha}(I_c)\subset I_c$ implies $\alpha(I)\subset I$.
Finally, it is obvious that $I\neq R$ (or $\{\Theta\}$) $\Leftrightarrow$ $I_c\neq A$ (or $\{\Theta\}$).
The proof is completed.
\end{proof}

Let us formulate one result from the work \cite{Rakhimov2018}.

\begin{pro}
  \label{pro:2}
A real C*-algebra $R$ is an injective if and only if the C*-algebra $R+iR$ is an injective.
\end{pro}

Now let's prove the following result.

\begin{pro}
  \label{pro:3}
Let $R$ be a real C*-algebra. Then a real C*-algebra $B$ is injective envelope of $R$
if and only if $B+iB$ is injective envelope C*-algebra of $A=R+iR$.
\end{pro}
\begin{proof}
Let $B+iB$ be an injective envelope of $R+iR$. By Proposition \ref{pro:2} \ $B$ is an injective real C*-algebra.
If $S$ is an injective real C*-algebra with $R\subset S$, then $A=R+iR\subset S+iS$ and by Proposition \ref{pro:2} $S+iS$ also is an injective, hence
$B+iB\subset S+iS$. Therefore $B\subset S$ and $B$ is injective envelope real C*-algebra of $R$.

Conversely, let $B$ be an injective envelope real C*-algebra of $R$. Then $A=R+iR\subset B+iB$ and by Proposition \ref{pro:2} \ $B+iB$ is an injective C*-algebra.
It's easy to see that $B+iB$ is injective envelope of $A$.
The proof is completed.
\end{proof}

Hence, using the result of \cite{Hamana79} we obtain the following corollary

\begin{cor}
  \label{cor:1}
Any real C*-algebra has a unique injective envelope real C*-algebra.
\end{cor}

\begin{pro}
  \label{pro:4}
Let $R$ be a real C*-algebra, $\alpha$ be a *-automorphism of $R$.
Then $\alpha$ is outer *-automorphism of $R$ if and only if $\tilde{\alpha}$ is outer *-automorphism of $A=R+iR$.
\end{pro}
\begin{proof}
If $\alpha$ is an inner *-automorphism of $R$, there is an unitary $u\in R$ such that $\alpha(x)=Adu(x)=uxu^*$, $\forall x \in R$.
Hence, we obtain
$$
\tilde{\alpha}(x+iy) = \alpha(x) + i\alpha(y) = uxu^* + iuyu^* = u(x+iy)u^* = Adu(x+iy),
$$
i.e., $\alpha$ is also an inner *-automorphism of $A$. Conversely, let $\tilde{\alpha}$ is an inner,
i.e. $\tilde{\alpha}(x+iy)=Adv(x+iy)$, where $v\in A$ is an unitary.
Since $\tilde{\alpha}(R)\subset R$, then by \cite[Corollary 3.1]{Kim2023} there  exists an unitary $u\in R$ such that $\tilde{\alpha}=Adu$.
Therefore $\alpha=Adu$, i.e., $\alpha$ is also an inner. The proof is completed.
\end{proof}

Now we will prove the main result of the section.

\begin{thm}
  \label{thm:1}
Let $R$ be real C*-algebra, $\alpha$ be a *-automorphism of $R$
such that $R$ is $\alpha$-simple.
Let $B$ be the injective envelope real C*-algebra of $R$. If $\alpha$ is an outer *-automorphism of $R$,
then $\alpha$ has a unique extension to an outer *-automorphism of $B$.
\end{thm}
\begin{proof}
By Proposition \ref{pro:1} C*-algebra $A=R+iR$ is $\tilde{\alpha}$-simple.
By Proposition \ref{pro:4} *-automorphism $\tilde{\alpha}$ also is an outer *-automorphism.
By Proposition \ref{pro:3} C*-algebra $B_c=B+iB$ is the injective envelope of $A$.
Then by \cite[Theorem 3.6]{Saito84} *-automorphism $\tilde{\alpha}$ has a unique
extension $\bar{\tilde{\alpha}}$ to outer *-automorphism of $B_c$.
It is obvious that the restriction of $\bar{\tilde{\alpha}}$ on $R$ coincides with $\alpha$, i.e.
$\bar{\tilde{\alpha}}|_R=\alpha$. Then it directly shows that $\bar{\alpha}=\bar{\tilde{\alpha}}|_{B}$ is a unique extension of $\alpha$ to an outer *-automorphism of $B$. The proof is completed.
\end{proof}

\section{Injective envelope of real simple C*-algebras}

To motivate the next definitions, suppose $A$ is a *-ring with unity, and let $w$ be a partial isometry in $A$.
If $e=w^*w$, it results from $w=ww^*w$ that $wy=0$ iff $ey=0$ iff $(1-e)y=y$ iff $y\in (1-e)A$, thus the elements
that right-annihilate $w$ form a principal right ideal generated by a projection.
If $S$ is a nonempty subset of $A$, we write
$R(S)=\{x\in A: sx=0, \ \forall s\in S\}$ and call $R(S)$ the {\it right-annihilator} of $S$.
Similarly, the set $L(S)=\{x\in A: xs=0, \ \forall s\in S\}$ denotes the {it left-annihilator} of $S$.
A {\it Baer *-ring} is a *-ring $A$ such that, for every nonempty subset $S$ of $A$, $R(S)=gA$
for a suitable projection $g$. It follows that $L(S)=(R(S^*))^*=(hA)^*=Ah$ for a suitable projection $h$.
A real (resp. complex) {\it AW*-algebra} is a real (resp. complex) C*-algebra that is a Baer *-ring
(for more details see \cite{Berberian2011}).
An AW*-algebra $A$ is called a {\it factor} if the center of $A$ is trivial.
It is known that, every W*-algebra is an AW*-algebra.
The converse of it was shown to be false by J.Dixmier, who showed that exist commutative AW*-algebras
that cannot be represented (*-isomorphically) as W*-algebras on any Hilbert space.\\[-2mm]

The following interesting result is true

\begin{thm}
  \label{thm:1}
Let $R$ be a real  C*-algebra with the unit and let $B$ its injective envelope.
If $R$ is a simple algebra, then $B$ is also simple, in which case $B$ is a real AW*-factor.\\[-2mm]
\end{thm}
\begin{proof}
By Proposition \ref{pro:3} C*-algebra $B+iB$ is an injective envelope of $A=R+iR$.
By Theorem 3.1, since $R$ is a simple algebra, then $A$ is also simple. By \cite[Proposition 4.15]{Hamana79}
$B+iB$ is also simple and $B+iB$ is a (complex) AW*-factor. Then $B$ is also simple and
by \cite[Proposition 4.3.1]{Ayupov2010} algebra $B$ is a real AW*-factor.
The proof is completed.
\end{proof}

Now, by analogy with \cite{Hamana79}, we give an example of an injective non W*-, AW*-factor of type III.\\[-2mm]

Consider the Calkin algebra: $A=B(H)/K(H)$, where $H$ is an infinite-dimensional separable Hilbert space
and $K(H)$ is an algebra of all compact operators on $H$.
In \cite[Example 5.1]{Hamana79} shows that if $B$ is an injective envelope of $A$, then $B$ is an injective AW*-factor of type III,
which is not a W*-algebra. Let us recall here \cite{Voiculescu76} that the Calkin algebra $A$ is not an AW*-algebra.

Now let's look at the real analogue of this example.
Let $H_r$ is an infinite-dimensional separable real Hilbert space.
Then since $B(H)=B(H_r)+iB(H_r)$ and $K(H)=K(H_r)+iK(H_r)$, where $H=H_r+iH_r$.
It's easy to see that $A_r=B(H_r)/K(H_r)$ is a real C*-algebra and we have $A_r+iA_r=A$.
Moreover, if $B_r$ is an injective real envelope of $A_r$, then $B_r+iB_r=B$.
By analogy \cite{Voiculescu76}, since $B$ is not a W*-algebra, then $B_r$ is also not real W*-algebra.

Thus $A_r$ is a real C*-algebra that is not real AW*-algebra and
the injective envelope $B_r$ of $A_r$ is an injective real AW*-factor of type III,
which is not a real W*-algebra.\\[-2mm]

It is known \cite{Connes76} that for each number $\lambda\in[0,1]$ the class of injective (complex) W*-factors of type III$_\lambda$
is unique, i.e. any two injective III$_\lambda$-factors are isomorphic.

In the real case: up to isomorphism there exist exactly two injective real W*-factors of type III$_\lambda$ ($0<\lambda\leq <$)
and up to isomorphism there exists a unique injective real W*-factor of type III$_1$.
For the case of real injective type III$_0$ factor we can construct a countable number of pairwise non isomorphic
real injective factors of type III$_0$, with isomorphic enveloping (complex) W*-factors (see \cite{Ayupov2010}, \cite{Ayupov97}).\\[-2mm]

Thus, from the above example and from the Hamana's example \cite[Example 5.1]{Hamana79} we get the following interesting result.

\begin{cor}
  \label{cor:2}
Up to isomorphism, the class of injective complex (real) AW*-factors of type III is at least one larger than the class
injective complex (respectively, real) W*-factors of type III.
\end{cor}

\end{document}